\documentclass[twocolumn,preprintnumbers,showpacs,pre,floatfix,amsmath,amssymb,notitlepage,superscriptaddress]{revtex4-1}

\usepackage{bm,bbm}
\usepackage{graphicx,epsfig,subfigure}
\usepackage{amsmath,amssymb}
\usepackage{amsthm}
\usepackage{color,xcolor}
\usepackage{url}
\usepackage[utf8]{inputenc}
\usepackage{enumitem}
\usepackage{setspace}
\usepackage{dcolumn}
\usepackage{bm}
\usepackage{times}
\usepackage{enumerate}
\usepackage{footnote}
\input epsf
\usepackage[titletoc]{appendix}

\usepackage{mathbbol}
\usepackage{bm}
\usepackage{bbm}
\usepackage{graphicx}
\usepackage{url}
\usepackage{xcolor}
\usepackage{enumitem}

\begin{document}

\newtheorem{thm}{Theorem}[section] 
\newtheorem{cor}{Corollary}[section] 
\newtheorem{definition}{Definition}[section]
\newtheorem{remark}[thm]{Remark}

\newcommand{\new}[1]{{\color{black}#1}}
\newcommand{\drt}[1]{{\color{black}#1}}


\title{Persistent Homology of Convection Cycles in Network Flows}

\author{Minh Quang Le}
\email{minhquan@buffalo.edu}
\affiliation{Department of Mathematics, University at Buffalo, State University of New York, Buffalo, NY 14260, USA}

\author{Dane Taylor}
\email{danet@buffalo.edu}
\affiliation{Department of Mathematics, University at Buffalo, State University of New York, Buffalo, NY 14260, USA}

\date{\today} 

\begin{abstract}
Convection is a well-studied topic in fluid dynamics, yet it is less understood in the context of  networks flows. 
Here, we incorporate techniques from topological data analysis (namely, persistent homology) to automate the detection and characterization of  convective/cyclic/chiral flows over networks, particularly those that arise for irreversible Markov chains (MCs).
As two applications, we study convection cycles arising under the PageRank algorithm, and we investigate chiral edges flows  for a stochastic model of a bi-monomer's  configuration dynamics.
Our experiments highlight how system parameters---e.g., the teleportation rate for PageRank and the transition rates of external and internal  state changes for a monomer---can act as \emph{homology regularizers} of convection, which we summarize with \emph{persistence barcodes} and  \emph{homological bifurcation diagrams.}
Our approach  establishes a new connection between the study of convection cycles and homology, the branch of mathematics that formally studies cycles, which has diverse potential applications throughout the sciences and engineering. 
\end{abstract}

\pacs{89.75.Hc, 02.50.Ga, 87.15.hj,84.35.+i}
\maketitle

\section{Introduction }

\drt{One of the main goals} of topological data analysis (TDA) is to characterize the structure of an object---usually a point cloud---through its   topological features. 
In particular, persistent homology \cite{edelsbrunner2010computational,otter2017roadmap} is a family of techniques that detect and summarize multiscale topological features and has been applied to a wide variety of  applications  including timeseries data \cite{khasawneh2016chatter,perea2015sliding}, 
image processing \cite{xin1994topology}, machine learning \cite{motta2019hyperparameter},  and
artificial intelligence \cite{gabella2019topology,liu2019scalable}.
%
Complementing the study of point-cloud data, another line of research involves utilizing the TDA toolset to study complex systems, for which   applications include
the analysis of spreading processes over social  networks \cite{taylor2015topological}, network neuroscience \cite{petri2014homological,giusti2015clique,kilic2022simplicial},  mechanical-force networks \cite{kondic2012topology}, jamming in granular material \cite{kramar2013persistence}, molecular structure \cite{liang1998analytical}, and DNA folding \cite{ichinomiya2020protein}.
%
%
In this paper, we employ techniques from TDA to 
study Markov chains (MCs), which provide a foundation to numerous areas of science and engineering including queuing theory \cite{kendall1953stochastic}, population dynamics \cite{kingman1969markov}, as well as statistical  (and machine learning) models that rely on Markov chain Monte Carlo \cite{gilks1995introducing}, hidden Markov models \cite{tierney1994markov}, and Markov decision process \cite{parr1998reinforcement}.

We utilize the mathematical framework of persistent homology to automate the detection  (and summarize the multiscale properties) of convection cycles that arise for the stationary flows of irreversible MCs.
Notably, while convection cycles have been extensively studied in fluid dynamics, they are less understood in the context of flows over networks. For example, it was recently observed that the coupling together of reversible MCs can give rise to an irreversible MC with  convection cycles that are an emergent property \cite{taylor2015topological}. Emergent convection cycles have also been recently found to describe the phenomenon of ``chiral edge flows'' \cite{tang2021topology}, providing new insights into the quantum Hall effect, configurational dynamics of monomers, and biological (e.g., circadian) rhythms. Given the inherent prevalence of convection cycles in MCs and other network flows, it is important that we place their study on a stronger mathematical, computational, and theoretical footing.





We study convection using a branch of mathematics called \emph{homology} and the related field   \emph{computational homology} \cite{kaczynski2006computational}. Both are concerned with studying the absence/presence of $k$-dimensional ``holes'' (and their connectivity) within a topological space such as a simplicial complex. Importantly, cycles on a graph are 1-dimensional (1D) holes, and so persistent homology is a natural fit to analyze convection cycles.
We construct filtrations of graphs by including edges according to the stationary flows across them (which is done in descending order so that the last edges to be included are those with the smallest stationary flows), and we summarize the persistent homology of the filtered graphs' associated clique complexes. See Fig.~\ref{fig:undirect_g} for a graph and its associated clique complex.
%
Computationally, we implement these techniques by building on a popular TDA framework called Gudhi \cite{gudhi:urm}, which we adapt to implement \emph{edge-value clique (EVC) filtrations} of scalar-functions that are defined over the the edges of a graph.

We apply this technique to two applications.
First, we study convection cycles arising under the PageRank algorithm \cite{page1999pagerank}, examining the role of  the \emph{teleportation parameter}. 
Second, we study chiral edge flows that emerge for a 4-state model that describes the configurational changes of a bi-monomer \cite{tang2021topology}, examining the roles of the external and internal transition rates.
These parameters significantly affect  convection  cycles arising for these respective applications, and we show that they act as  ``homology regularizers'' of convection. We introduce ``homological bifurcation diagrams'' to summarize these effects. Our methods provide mathematically principled (and automated) tools to gain a deeper understanding of the structural patterns of convection on networks, and they are expected to be  useful to myriad applications across the physical, social, biological and computational sciences.


The remainder of this paper is organized as follows:
We present background information in Sec.~\ref{sec:back}, 
\drt{our methodology} in Sec.~\ref{sec:results1}, 
\drt{applications in Sec.~\ref{sec:apps},}
and a discussion in Sec.~\ref{sec:conclusion}.

\section{Background Information }\label{sec:back}

Here, we  present introductory material about 
simplicial complexes and homology (Sec.~\ref{sec:simplex}), 
persistent homology of graphs (Sec.~\ref{sec:PHw}), 
and discrete-time MCs (Sec.~\ref{sec:conv_cycle}).

\subsection{Simplicial complexes (SCs) and their homology}\label{sec:simplex}

We first define an \emph{undirected graph} $G = \{\mathcal{V},\mathcal{E}\}$, where $\mathcal{V}=\{1,\dots,N_0\}$ is a  set of $N_0$ vertices and $\mathcal{E}\subset \mathcal{V}\times \mathcal{V}$ is a set of edges. 
Note that each vertex is specified by a single index $i\in V$, and each  edge is specified by an unordered pair  $(i,j)\in \mathcal{V}\times \mathcal{V}$. 
More generally, we define a $(k+1)$-tuple of vertices $\sigma = (i_0,i_1,\dots,i_{k})$ as a {$k$-dimensional simplex}, or $k$-\emph{simplex}
\footnote{The notion of ``dimension'' is more interpretable for the case of a non-abstract simplicial complex, for which the vertices correspond to locations in a Euclidean metric space. In that case, every $k$-simplex is defined as the $k$-dimensional surface that is contained by its faces, which themselves are $(k-1)$ dimensional surfaces. For example, a 2-simplex is a triangle defined as the interior of 3 line segments (which are the 3 1-dimensional cofaces of the 2-simplex).}. Vertices  and edges are equivalent to 0-simplices and 1-simplices, respectively. 
%
An abstract SC is a set of simplices of possibly different dimensions, and it is a generalization of an undirected graph. It is also a type of \emph{hypergraph} with a constraint on which simplices can exist.
That is, for any $k$-simplex $(i_0,i_1,\dots,i_{k})$ in an SC, its \emph{faces} are the $(k-1)$ simplices in which one of the of indices is omitted (e.g., $i_1$ is omitted to yield $(i_0,i_2,\dots,i_{k})$). The \emph{cofaces} of a $(k-1)$-simplex are the $k$-simplices for which it is a face.
%
%
Note that the faces of an edge $(i,j)$ are the vertices $i$ and $j$, and likewise, $(i,j)$ is a coface of each of these vertices.

With these definitions, we state the two  restrictions that are required for a SC:
(i) for any face, its faces must  be included in the SC;
and (ii) the intersection of any two faces is either a face of both, or it is an empty set. 
The dimension of an SC is the maximum dimension of its simplices, and an undirected graph is a 1-dimensional SC---it contains 0-simplices and 1-simplices, and for any edge $(i,j)$, the vertices $i$ and $j$ must exist.
We will focus on a particular type of SC that can be generated from a graph and is called a \emph{clique complex}.
A clique complex $K(G)$ of a graph $G$ is the SC in which there is a 1-to-1 correspondence between the $(k+1)$-cliques in the graph and the $k$-simplices in the SC. (Recall that an \emph{n-clique} is a complete subgraph on $n+1$ vertices of a graph.) Given this 1-to-1 correspondence, the map from $G$ to $K(G)$ is invertible, and $G$ can be recovered as the $1$-skeleton of $K(G)$. (A \emph{$k$-skeleton} of a SC is the SC that is obtained after removing all simplices  having dimensions that are greater than $k$.)
%

%
%
%
%

\begin{figure}[t]
\centering
\includegraphics[width=\linewidth]{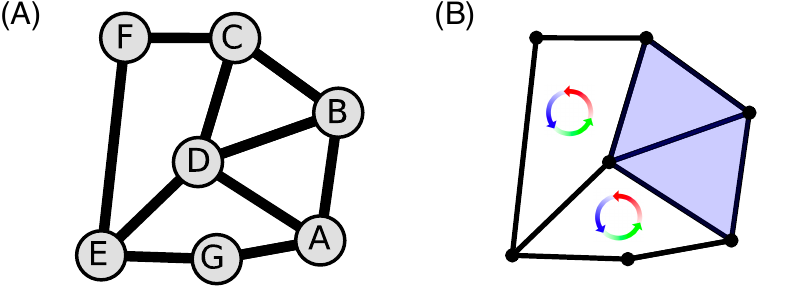}
\caption{{\bf A graph and an associated simplicial complex (SC).}
(A)~Graph $\mathcal{G}$ with $N=7$ vertices and $M= \drt{10}$ undirected, unweighted edges.
(B)~The corresponding \emph{clique complex} $\mathcal{S}$, where each $k$-clique gives rise to a $(k-1)$-simplex. Each triangle (i.e., $3$-clique) gives rise to a $2$-simplex (see shaded triangles). 
The SC is a topological space and has  associated vector spaces.
%
Consider the space $\mathbb{R}^7$ of real-valued functions defined over the vertices. Since there is {just 1} connected component, the SC's \drt{0-dimensional (0D)} homology is a {1D subspace} of $\mathbb{R}^7$.
Similarly, there are \drt{two ``homological'' 1-cycles} that are not a boundary of a 2-simplex, and so the 1D homology is a 2D subspace of $\mathbb{R}^{10}$ (i.e., the space of real-valued functions defined over the ten edges).
%
}
\label{fig:undirect_g}
\end{figure}

\begin{figure*}[t]
\centering
\includegraphics[width=1\linewidth]{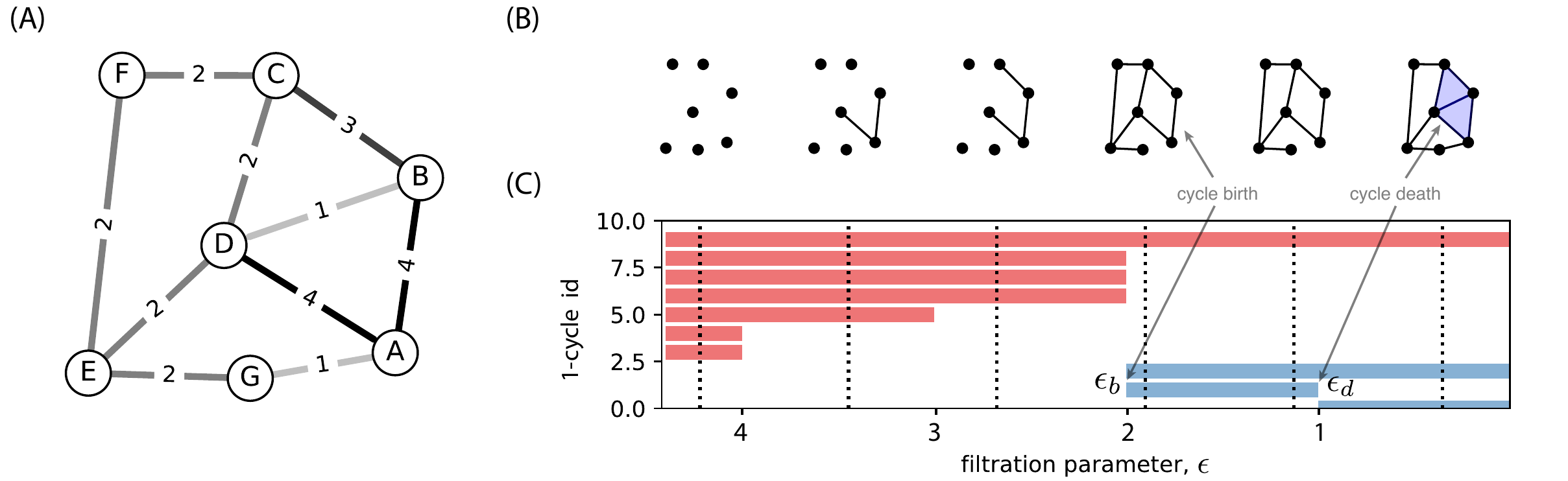}
\caption{
{\bf Persistent homology of a graph according to an \emph{edge-value clique (EVC) filtration}.}
(A) An undirected graph ${G}(\mathcal{V},\mathcal{E})$ with a scalar function $f:\mathcal{E}\to \mathbb{R}$ defined over the edges $\mathcal{E}$. 
We apply an EVC filtration to the graph  by considering a monotonically decreasing \emph{filtration  parameter} $\epsilon = (0,5]$, and by considering a filtered sequence of graphs ${G}(\mathcal{V},\mathcal{E}_\epsilon)$, where $\mathcal{E}_\epsilon = \{(i,j)|f(i,j)> \epsilon\}$ is the subset of edges for which $f(i,j)$   are larger than a threshold $\epsilon$.
(B)~Visualization of the graphs' associated clique complexes   $K_\epsilon \equiv K({G}(\mathcal{V},\mathcal{E}_\epsilon))$  for several $\epsilon$. 
%
%
%
%
%
(C)~A \emph{persistence barcode} summarizes how the 0-dimensional (red) and 1-dimensional (blue) 
\drt{homological $k$-cycles} of
$K_\epsilon$ change with decreasing $\epsilon$. 
The arrows highlight two events: 
at $\epsilon=\epsilon_b$, a \drt{homological} 1-cycle 
involving four edges is \emph{born}; at $\epsilon=\epsilon_d$, the \drt{1-cycle} \emph{dies}, since it is ``filled in'' by a 1-simplex and two 2-simplices.
}
\label{fig:filtra_g}
\end{figure*}

We next discuss \emph{simplicial homology}, \drt{which will lead to a formal definition of ``homological'' cycles. To this end, we} 
consider  vector spaces defined over the $k$-simplices in a SC.
%
A \emph{$k$-chain}, $\sum_{n=1}^{N_k} \alpha_n \sigma_n$, is a linear combination of $k$-simplices $\{\sigma_n\}$ with weights $\{\alpha_k\}$. (Note that $N_0$ and $N_1$ are the numbers of vertices and edges, respectively.) If a  SC contains $N_k$ different $k$-simplices, then the vector space of $k$-chains is $N_k$-dimensional, and it is \drt{isomorphic} to $\mathbb{R}^{N_k}$ if one assumes $\alpha_k\in\mathbb{R}$. %
We now consider a simplicial map $f: X_k\to X_{k-1}$ between $X_k$, which is a SC of dimension $k$, and $X_{k-1}$, which is a SC of dimension $k-1$ that contains the faces of simplices in $X_k$. 
%
%
%
Considering the vector space $C_k$ of $k$-chains defined over $k$-simplices in $X_k$ and vector space $C_{k-1}$ of $(k-1)$-chains defined over their cofaces in $X_{k-1}$, we define 
the linear \emph{boundary map} $\partial_k: C_k \longrightarrow C_{k-1}$, 
where the action of $\partial_k$ on any  $k$-simplex  is given by
\begin{align}
\partial_k(i_0,...,i_k)=\sum_{j=0}^k (-1)^j(i_0,...,i_{j-1},i_{j+1},...,i_k).
\end{align}
The boundary map allows one to relate vectors in $C_{k}$ to those in $C_{k-1}$.
For example, the \emph{boundary} of a 2-simplex (i.e., triangle) $(i,j,k)$ is the signed combination if the associated edges, $\partial_2 (i,j,k) = (j,k) - (i,k) + (i,j)$.
Notably, the boundary of any closed \drt{path} is  zero, 
\drt{which yields an algebraic definition of a  $k$-cycle: any $k$-chain that lies within the subspace $Z_{k}$, where}
$Z_{k}=\text{ker}(\partial_{k}) \subseteq C_k$ \drt{is}
the \emph{vector space  of  $k$-cycles}.

\drt{Notably, $k$-cycles can arise for different reasons, and we distinguish two types.}
The boundary map satisfies the property $\partial_{k}\circ \partial_{k+1} =0$, which essentially states that the boundary of a boundary is zero. [For the triangle, $\partial_1 \circ \partial_2 (i,j,k) = \partial_1(j,k) - \partial_1(i,k) + \partial_1(i,j) = (k - j) - (k-i)  + (j-i)=0 $.]
Thus we define $B_{k}=\text{image}(\partial_{k+1})$ as the \emph{subspace of $(k+1)$-boundaries}, and it follows that $B_{k}\subseteq Z_{k}$. In other words, some cycles arise simply because they are boundaries of $(k+1)$-simplices For example, observe in Fig.~\ref{fig:undirect_g}(B) that there are two ``triangular'' cycles that exist around the two 2-simplices, but that there are two other cycles that also exist. The \emph{$k$-th simplicial homology} is  defined as the quotient space $H_{k}= Z_{k}/B_{k}$, and it represents the subspace of \drt{$k$-dimensional   cycles (i.e., $k$-cycles)  that do not arise simply as the boundary of a $(k+1)$-simplex. } 

%

\drt{The $k$-th simplicical homology can be represented by the span of}
\emph{homology generators}, \drt{which are a linearly independent} set of $k$-chains that span $H_n$ and represent the associated $k$-cycles. The number of \drt{linearly independent} homology generators 
is called a \emph{Betti number}
\begin{align}
\beta_k = \dim H_k = \dim (Z_k) - \dim (B_k).
\end{align}
Informally, $\beta_0$ is the number of connected components;
$\beta_1$ is the number of 1-dimensional cycles or ``loops'' (that is, not including the triangular boundaries of 2-simplices); \drt{and}
$\beta_2$ is the number of 2-dimensional holes or ``voids'' (e.g., the interior of a triangulated sphere). 
For the SC shown in Fig.~\ref{fig:undirect_g}(B), $\beta_0=1$ since there's  one connected component, and $\beta_1=2$ since there are two cycles that are  not  simply  the boundaries of  $2$-simplices.

\drt{
By formulating $k$-cycles algebraically, one can consider the linear dependence and independence of $k$-cycles. As such, one can not only identify cycles, but also investigate the relations/connectivity between cycles,  which we find to be instrumental for understanding  pattern formation for cycle. 
We also highlight that a given homological $k$-cycle can potentially have more than one homological generator. Such generators are said to be \emph{homologically equivalent}, and they can be obtained by considering linear combinations of $k$-cycles (including both homological $k$-cycles and boundaries). We will later show that this complicates the investigation of  convection cycles through the lens of homological $k$-cycles.
}





%
%

\subsection{Persistent homology of scalar functions defined over edges}\label{sec:PHw}

One of the greatest tools of topological data analysis is the study of \emph{persistent homology} \cite{edelsbrunner2010computational,otter2017roadmap}. Here, we examine how the homology of a topological object  changes as it undergoes a \emph{filtration} to  yield a monotonically increasing sequence  $X_0 \subseteq X_1 \subseteq X_1 \subseteq ...$ (e.g., of simplicial complexes $\{X_t\}$).
We consider filtrations in which one has a {scalar  function} $f:\mathcal{E} \longrightarrow \mathbb{R} $ over the edges, and each edge $(i,j)\in\mathcal{E}$ is retained/removed according to $f(i,j)$. The values $f(i,j)$ could be edge weights for a weighted graph, but in general they  can  encode any scalar property. We visualize such a graph and the values $f(i,j)$ in Fig.~\ref{fig:filtra_g}(A).

We call the process  an \emph{edge-value clique (EVC) filtration}, and we construct it as follows. Given a graph $G(\mathcal{V},\mathcal{E})$ and a \emph{filtration function} $f $, we   define   the subsets $\mathcal{E}_\epsilon = \{(i,j)|f(i,j)> \epsilon\}$    in which one retains  edges only for which $f(i,j)$ is sufficiently large. Note that the subsets $\{\mathcal{E}_\epsilon \}$ are non-decreasing  as $\epsilon$ decreases (i.e., $\mathcal{E}_\epsilon \subseteq \mathcal{E}_{\epsilon'} $ for any $\epsilon'<\epsilon$). We must specify a range over which to decrease $\epsilon \in(\epsilon_B,\epsilon_A]$, and in practice we assume
$\epsilon_A >\max_{(i,j) \in \mathcal{E}} f(i,j)$ and $\epsilon_B<\min_{(i,j) \in \mathcal{E}} f(i,j)$. It then follows that $\mathcal{E}_{\epsilon}=\emptyset$ is an empty set of edges when $\epsilon \ge \epsilon_A$, and $ \mathcal{E}_\epsilon =  \mathcal{E}$ (i.e., all edges are retained) when $\epsilon \le \epsilon_B$.
See \cite{minh_code} for our codebase that implements EVC filtrations by adapting the     TDA framework called Gudhi \cite{gudhi:urm},  and which reproduces the results of this paper.

\begin{figure*}[t]
\centering
\includegraphics[width=\linewidth]{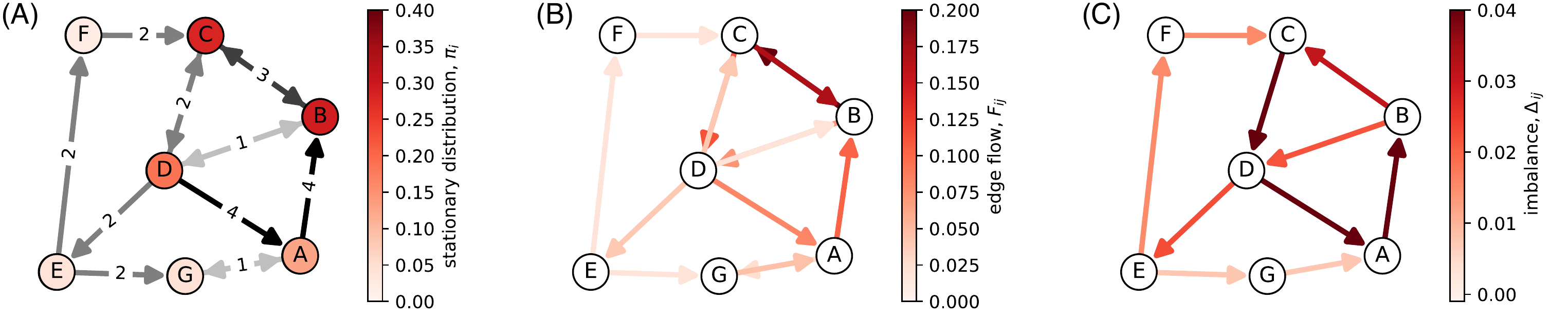}
\caption{
{\bf Stationary distribution, edge flows, and convection cycles for an irreversible MC.}
%
We study a discrete-time random walk over a directed, weighted graph that resembles the undirected graph  in Fig.~\ref{fig:filtra_g}, except that the edges are now either directed or bidirectional. (Recall that undirected graphs give rise to reversible MCs that lack convection cycles.)
(A) The color of each vertex indicates the \emph{stationary distribution} $\pi_i$ of random walkers at each vertex $i$.
(B) Edge colors indicate the \emph{stationary  flows} $F_{ij}=\pi_iP_{ij}$ across edges, i.e., the stationary fraction of random walkers that traverse each directed edge.
(C)  
\emph{Flow imbalances} $\Delta_{ij} =   (F_{ij} - F_{ji})$ manifest as a pattern of \emph{convection cycles}.
By construction, $\Delta_{ij}=-\Delta_{ji}$, and we use arrows to indicate the directions of imbalance, e.g.,
$i\rightarrow j$ if $\Delta_{ij}>0$.  
 }
\label{fig:direct}
\end{figure*}

In Fig.~\ref{fig:filtra_g}(B), we visualize a sequence of filtered clique complexes $\{ K_\epsilon\}$   that are associated with the  filtered graphs $\{G_\epsilon\}$ that are defined with the  edge  sets $\{\mathcal{E}_\epsilon\}$.
In Fig.~\ref{fig:filtra_g}(C), we summarize the persistent homology of $\{ K_\epsilon\}$  in  a \emph{persistence barcode}, which   reveals how  homology changes with $\epsilon$.
Observe that when $\epsilon $ is sufficiently large,   $K_\epsilon$ contains vertices but no edges. On the other hand, when $\epsilon$ decreases to be sufficiently small, then $K_\epsilon$ recovers the original clique complex [recall Fig.~\ref{fig:undirect_g}(B)]. The values of $\epsilon$ that were used to create Fig.~\ref{fig:filtra_g}(B) are indicated by the vertical dotted lines in Fig.~\ref{fig:filtra_g}(C).

Each \drt{horizontal} bar in the persistence barcode \drt{shown   in Fig.~\ref{fig:filtra_g}(C)} indicates the \emph{lifetime} of a 
\drt{homological 1-cycle---that is,}  the values of $\epsilon$ for which it exists. The red and blue bars reflect $0$-homology and $1$-homology respectively. The dimensions of the homology spaces (i.e., Betti numbers) can be found by counting the number of 
\drt{homological 1-cycles}
at any particular $\epsilon$. For example, one can observe that $\beta_1=0$ when $\epsilon=3.5$, $\beta_1=1$   when $\epsilon=2.5$, and   $\beta_1=2$ when $\epsilon=1.5$. Clearly, the \drt{homological $1$-cycles}
are undergoing bifurcations as $\epsilon$ varies.
A persistence barcode is convenient to identify for each generator:
the value $\epsilon_b$ of $\epsilon$ when it is ``born'' (i.e., the \drt{homological $k$-cycle}
does not exist when $\epsilon>\epsilon_b$);
the value $\epsilon_d$ of $\epsilon$ when it  ``dies'' (i.e., the \drt{homological $k$-cycle}
does not exist when $\epsilon<\epsilon_d$);
its \emph{lifetime}  $ (\epsilon_d,\epsilon_b]$;
and   \emph{lifespan} $|\epsilon_d-\epsilon_b|$.
A cycle's {lifespan}  quantifies its persistence under the filtration, and it is  often interpreted as a proxy for the cycle's significance
\drt{(although  short-lifetime cycles can also be   important in certain contexts).}
\subsection{Convection cycles for irreversible Markov chains (MCs)}\label{sec:conv_cycle}

We will apply persistence homology to study convection cycles in irreversible MCs \cite{lovasz1993random}, which we now briefly summarize.
%
%
A discrete-time MC is a ``memoryless'' random process in which for time steps $t=0,1,2,\dots$ ,
the system state $S_{t}\in\mathcal{V}$ satisfies the Markov property 
${ \text{P}}[S_{t+1}=i |S_0=i_0, ..., S_t=i_t]
= {  \text{P}}[S_{t+1}=i|S_t=i_t]$, which implies that the probability of a state occurring at the next time step only depends on the current state and not earlier states.
In our case, we consider MCs that correspond to a random walk on a (possibly) weighted and  directed graph having an \emph{adjacency matrix}  ${\bf A}$ in which $A_{ij}\in\mathbb{R}$ is nonzero if $(i,j)$ is an edge, $(i,j)\in\mathcal{E}$, and $A_{ij}=0$ otherwise. 
We similarly define a \emph{transition matrix},  ${\bf P} = {\bf  D}^{-1} {\bf A}$, where ${\bf D}$ is a diagonal matrix  with entries that encode the  (possibly) weighted vertex degrees ${\bf D}_{ii} = \sum_j A_{ij}$. For directed graphs, each $(i,j)$ is considered to be an ordered pair, and each ${\bf D}_{ii}$ encodes the  out-degree of vertex $i$.
Each matrix element $P_{ij}$ gives the probability for a random walk to transition from   vertex  $i$ to $j$.
%
Letting $x_i(t)$ denote the probability that the system is in state $i$  (or equivalently, the probability that a random walker is at vertex $i$) at time $t$, one can utilize the Markov property to obtain the linear discrete-time system $x_j(t+1) = \sum_i x_i(t)P_{ij}$. By defining ${\bf x}(t) =[x_1(t),\dots,x_{N_0}]^T$,   one equivalently has 
\begin{align}
    {\bf x}(t+1)^T = {\bf x}(t)^T {\bf P}.
\end{align} 
Since ${\bf x}(t)$ is a vector of probabilities, we assume that it is normalized in 1-norm, $\sum_i x_i(t)=1$.

Herein, we focus on network flows after a system reaches a \emph{stationary state}, in which case ${\bf x}(t)$ converges to a limiting vector $\pi = \lim_{t\to\infty}{\bf x}(t)$ that satisfies the eigenvalue equation $\pi^T = \pi^T {\bf P}$. By construction, $\pi$ is a vector of probabilities and contains nonnegative entries. Furthermore, as a row-stochastic matrix, ${\bf P}$ has an eigenvalue equal to one (i.e., the largest eigenvalue) and its right dominant eigenvector is the vector containing 1's as entries. Our assumption of convergence requires that matrix ${\bf P}$ is an irreducible and 
aperiodic   \cite{bapat1998max}  or that the initial condition ${\bf x}(0)$ lies in  a converging subspace. In the stationary state, the \emph{stationary  flow} across each edge $(i,j)$ per time step is given by
\begin{align}
    F_{ij} = \pi_i P_{ij}.
\end{align}

\begin{figure*}[t]
\centering
\includegraphics[width=\linewidth]{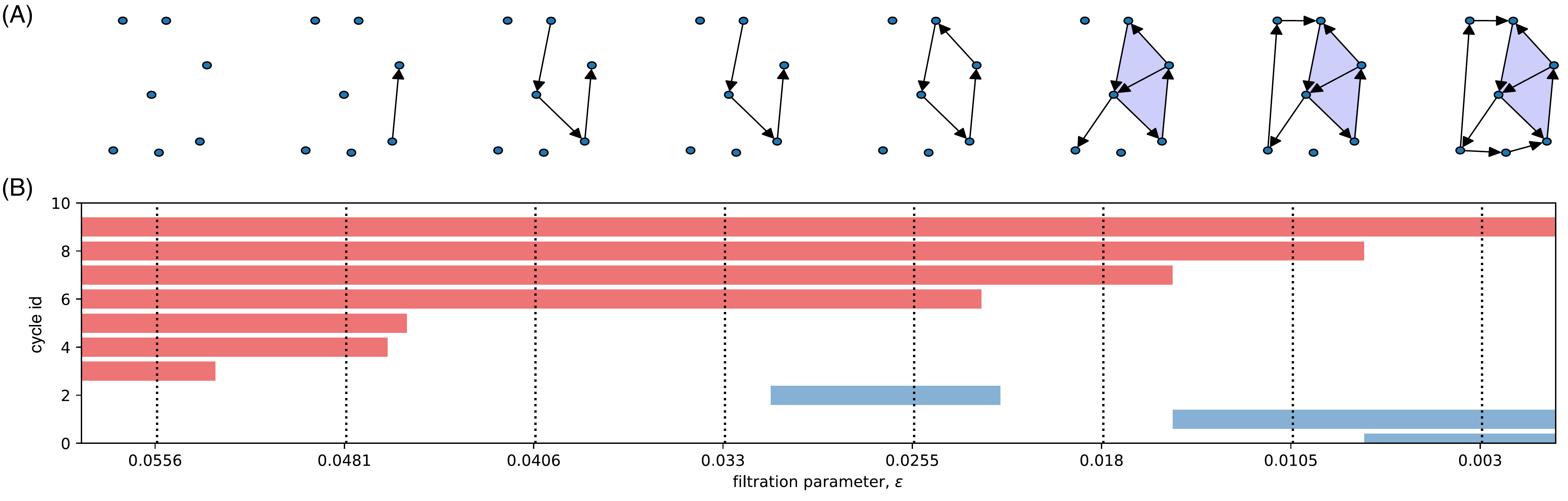}
\caption{{\bf Persistent homology of convection cycles.}
\drt{(A) Visualization of an EVC filtration applied to flow imbalances arising for} the irreversible MC shown in  Fig.~\ref{fig:direct}, and we use the magnitude $|\Delta_{ij}|$ of flow imbalance as the filtration function $f:\mathcal{E}\to\mathbb{R}$. \drt{We indicate flow imbalances' directions with arrows, noting that the clique complexes that are constructed by the filtration are undirected, since the filtration does not incorporate information about edge directions.}
\drt{(B) Persistence barcodes for homological 1-cycles.}
Observe that  the 1-cycle that first appears dies before the other 1-cycles are born.
}
\label{fig:Pbnet}
\end{figure*}

We study \emph{convection cycles} using an approach that was developed in \cite{taylor2020multiplex}.
Specifically, for each edge we define the stationary \emph{flow imbalance}
\begin{align}
\Delta_{ij} = F_{ij} - F_{ji}.
\end{align}
By construction, $\Delta_{ji} = -\Delta_{ij}$, and we say that the \emph{imbalance direction} is from $i$ to $j$ when $\Delta_{ij}>0$.
Importantly, the defining feature of a \emph{reversible MC} is that $\Delta_{ij} =0$ for all $i$ and $j$. That is, the directional flows match $\pi_i P_{ij}=\pi_j P_{ji}$ for any edge $(i,j)$. This is the case for any undirected graph, since in this case ${\bf A}={\bf A}^T$, and it follows that  $\pi_i = D_{ii}/\sum_j D_{jj}$. 
%
%
%
In contrast, an \emph{irreversible MC} yields asymmetric stationary flows and $\Delta_{ij}$ is nonzero for some edges. 
%
\drt{
To formally define convection cycles, we consider a new  graph $G_{\Delta}(\mathcal{V},\mathcal{E}_\Delta)$ such that  each positive value $\Delta_{ij}$ gives rise to a directed edge $(i,j,\Delta_{ij})\in\mathcal{E}_\Delta$ having weight $\Delta_{ij}$. We then define a \emph{convection cycle} to be any non-intersecting closed path in $G_{\Delta}(\mathcal{V},\mathcal{E}_\Delta)$.
}

\drt{
In Fig.~\ref{fig:direct}, we  illustrate for  an example MC how flow imbalances  manifest as a pattern of convection cycles. 
In Figs.~\ref{fig:direct}(A), 
\ref{fig:direct}(B), and 
Fig.~\ref{fig:direct}(C), we use edge colors to indicate the stationary distribution $\pi$, stationary edge flows $F_{ij}$, and flow imbalances $\Delta_{ij}$, respectively. 
Observe that some of the arrows in Figs.~\ref{fig:direct}(A)--(B) are bidirectional, since some of the graph's edges are bidirectional. In contrast, the arrows in Fig.~\ref{fig:direct}(C) are exclusively directed since they now indicate the directions of flow imbalances. There exists an edge $i\to j$ only if  $\Delta_{ij}>0$, which also implies $j\to i$ is not an edge since  $\Delta_{ji}=-\Delta_{ij}$. Observe in Fig.~\ref{fig:direct}(C) that this yields five convection cycles. In Sec.~\ref{sec:interp}, we will further discuss  these convection cycles and  their relation to homological $1$-cycles.}

\drt{
Before continuing, we highlight that convection cycles revealed through flow imbalances \cite{taylor2020multiplex} do not take into account the probability of transitioning to or away from a convection cycle, and so they are 
not necessarily ``cyclic traps.'' That is, the presence of a convection cycle does not imply that it is unlikely for a random walker to leave (or move in an opposite direction as) the cycle.
For example, observe in Fig.~\ref{fig:direct} that the counter-clockwise flow around convection cycle $A\to B\to C\to D$ is approximately 0.025, yet there is a flow of approximately $0.02$ that leaves the cycle at node D, and a flow of approximately $0.15$ moves in the opposite direction from node C to B. 
%
%
Future research will likely uncover complementary notions of convection with different advantages/disadvantages, and  our proposed techniques using persistent homology  can likely be similarly  extended.
}

\section{Homological analyses of convection  }\label{sec:results1}

We now employ  persistent homology to automate the detection, characterization, and summarization of the  homological patterns of  convection cycles. 
In Sec.~\ref{sec:r1_a}, we study the MC that was presented in   Fig.~\ref{fig:direct}.
\drt{In Sec.~\ref{sec:interp}, we discuss the relation between convection cycles and homological $1$-cycles.}

\begin{figure*}[t]
\centering
\includegraphics[width=\linewidth]{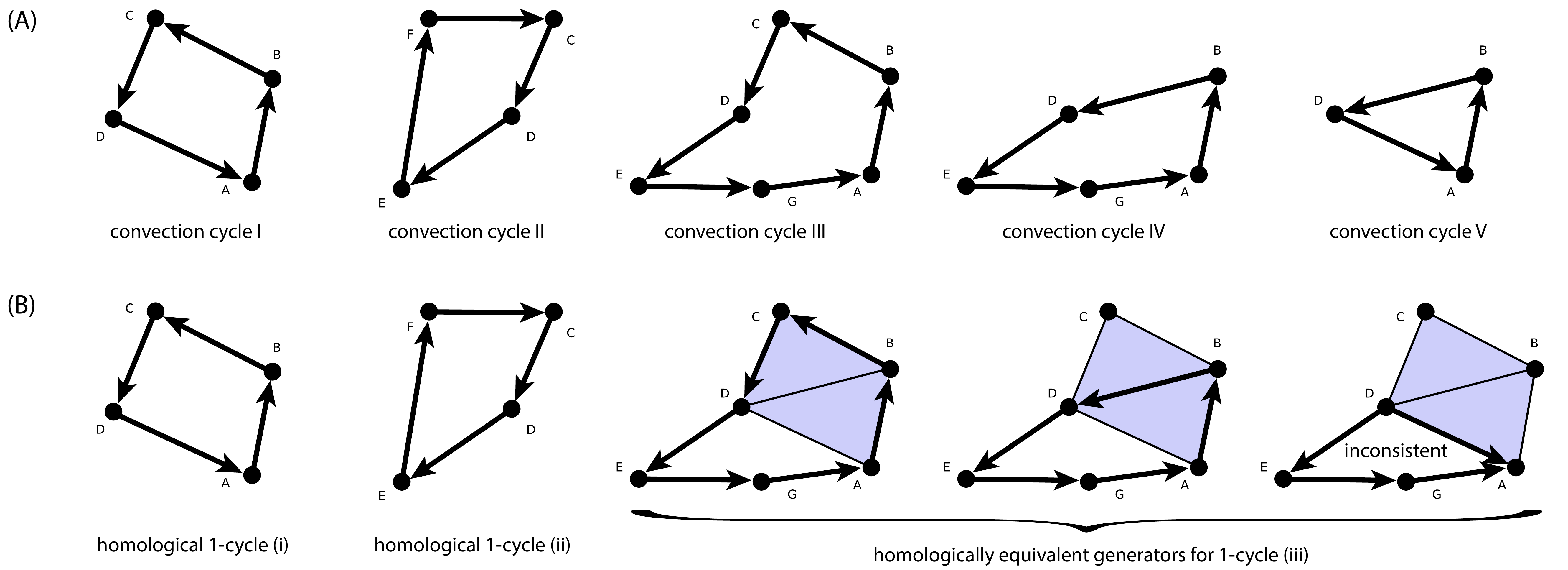}
\caption{
\drt{{\bf Relation between convection cycles and homological 1-cycles.} 
(A)~The flow imbalances  shown in Fig.~\ref{fig:direct}(C) give rise to five convection cycles, which we label I--V.
(B)~Persistent homology using EVC filtrations applied to a  network of flow imbalances reveals three \emph{homological 1-cycles}, which we label (i)--(iii). Each homological 1-cycle represents a ``1-dimensional hole'' and can be represented by one or more \emph{homological generator} (recall Sec.~\ref{sec:simplex}).
Observe that there is a one-to-one correspondence between convection cycles I and II and homological 1-cycles (i) and (ii). In contrast,  there are three homologically equivalent generators for 1-cycle (iii) as shown. Two of the generators correespond to convection cycles III and IV. The third generator does not correspond to a convection cycle, because the edge directions are not consistently in the same orientation (i.e., always clockwise or counter-clockwise).
%
}}
\label{fig:compare}
\end{figure*}

\subsection{Persistent homology of convection cycles}\label{sec:r1_a}

Recall from Sec.~\ref{sec:PHw} that EVC filtrations were defined for an undirected graph with a scalar function defined on the edges. Therefore,  \drt{given} an MC corresponding to a (potentially) directed and weighted graph, we first consider the associated undirected graph. Then we study homology under an EVC filtration in which the filtration function $f:\mathcal{E}\longrightarrow \mathbb{R}$ is given by the magnitudes of the  flow imbalances
\begin{align}
    f(i,j) = |\Delta_{ij}|.
\end{align}
In this way, the persistent homology that is revealed corresponds to the convection cycles that arise under flow imbalances.



In Fig.~\ref{fig:Pbnet}, we visualize persistence barcodes for an EVC filtration associated with the convection cycles shown in Fig.~\ref{fig:direct}(C). 
Note that this figure is analogous to Fig.~\ref{fig:filtra_g}(C), where we had previously chosen the filtration function $f(i,j)$ to be the edge weights. 
Since we now use a different function $f$,   the cycles  now have different births, deaths, lifetimes and lifespans.
Interestingly,  the 1-cycle involving vertices $\{A,B,C,D\}$ is now born and dies before the other two 1-cycles are  born. While there is an obvious connection between the EVC homology of a graph induced by edge weights and that which is induced by convective flows, this relation remains unclear and should be explored in future work.

\drt{We note that one could also construct EVC filtrations by increasing $\epsilon$ and retaining  edges $(i,j)$ for which $|\Delta_{ij}|$ is smaller than $\epsilon$. In Appendix~\ref{sec:appendix}, we provide an example illustrating why EVC filtrations with decreasing $\epsilon$ are superior to those with increasing $\epsilon$ for the goal of studying convection cycles. In particular, EVC filtrations that   decrease $\epsilon$ focus on  1-cycles that are associated with large-flow convection cycles (i.e., large values of $|\Delta_{ij}|$), which we consider to be the ones that are more significant. In contrast, EVC filtrations that  increase $\epsilon$ focus on  1-cycles that are associated with small-flow convection cycles (i.e., small values of $|\Delta_{ij}|$), which we consider to be less significant. 
}

\subsection{Comparing convection cycles and homological 1-cycles}\label{sec:interp}

\drt{
We propose to study pattern formation for convection cycles using persistent homology techniques for homological 1-cycles; however, one should keep in mind that these are two different notions for cycles.
Homological $1$-cycles are 1-dimensional holes for a topological space, and $k$-cycles generalize to higher dimensional by representing higher-dimensional holes (Sec.~\ref{sec:simplex}). 
In contrast, we define convection cycles to be closed non-backtracking paths in a directed graph that encodes flow imbalances (Sec.~\ref{sec:conv_cycle}).
In this section, we will clarify the relationship between homological $k$-cycles and convection cycles, thereby revealing the capabilities and disadvantages of   existing persistent homology techniques for  studying convection cycles.
Continuing with the previous example [see Figs.~\ref{fig:direct}-\ref{fig:Pbnet}], we find that  flow imbalances give rise to five convection cycles, which we enumerate I--V and visualize in  Fig.~\ref{fig:compare}(A).
In contrast, we identify three homological 1-cycles using persistent homology with EVC filtrations, which we enumerate (i)--(iii) and visualize in  Fig.~\ref{fig:compare}(B). 

Observe that there is a one-to-one correspondence between convection cycles I and II and homological 1-cycles (i) and (ii). 
Also observe that homological 1-cycle (iii) has three homologically equivalent generators, and any of them can be used to represent the 1-cycle (which again, is  defined as a 1-dimensional hole). Each subsequent generator can be obtained via a \emph{topological retraction} in which a 2-simplex is collapsed down onto one of its edges. Interestingly, the first two homological generators for 1-cycle (iii) correspond to convection cycles III and IV. In contrast,  the third  generator corresponds to a loop that is not a convection cycle, since the flow-imbalances' directions do not point in a consistent direction along the cycle (i.e., clockwise or counter-clockwise). 
Finally, observe that convection cycle V is a boundary of a 2-simplex, and it therefore does not contribute to the 1-dimensional simplicial homology. 

Thus, it is important to not misinterpret one notion of cycle for the other. At the same time, our findings in Fig.~\ref{fig:compare} also highlight that there is a need for new persistent homology techniques that cater specifically to convection cycles and directed graphs. For example, if one were to omit  the 2-simplex that involves vertices A, B and D from the clique complexes that arise  under an EVC filtration, then  convection cycle V would coincide with a homological 1-cycle.
However, the aim of this paper is not to develop new methods for persistent homology. 
Instead, we proposed to begin this pursuit by   studying convection cycles using existing methods for persistent homology.
Even   though there is not an exact one-to-one match between convection cycles and homological 1-cycles, because they are closely related, we find that persistent homology can effectively   detect and summarize  convection cycles' patterns.
}

\section{Applications}\label{sec:apps}

\drt{
In this section,  we apply our approach to two applications.}
In Sec.~\ref{sec:r1_b}, we study MCs arising for the Google PageRank algorithm, exploring how convection cycles are effected by the \drt{teleportation} parameter $\alpha$.  
In Sec.~\ref{sec:r1_c}, we study a type of emergent convection cycle called a {chiral edge flow}.

\begin{figure*}[t]
\centering
\includegraphics[width=1\linewidth]{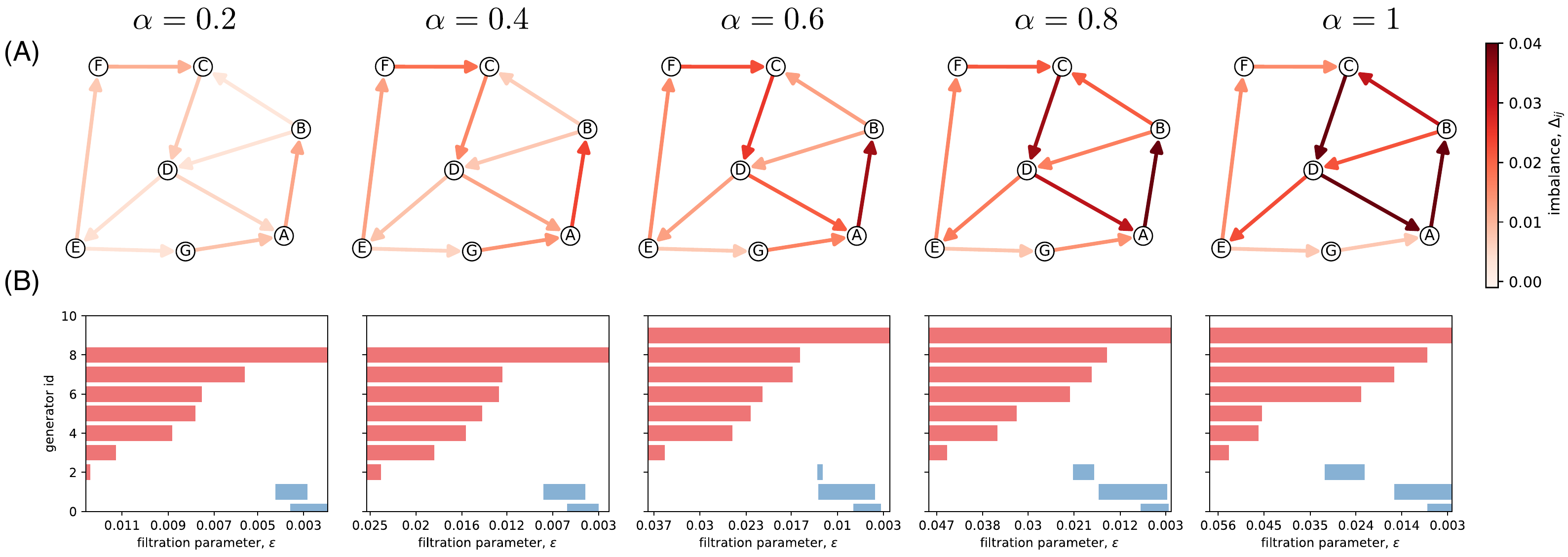}
\caption{
{\bf Persistent homology of convection cycles arising under PageRank.} 
We study the MC associated with PageRank for the directed graph from Fig.~\ref{fig:direct}(A) under  several choices of $\alpha$.
(A) Flow imbalances $\Delta_{ij}(\alpha)$ give rise to convection cycles. \drt{For clarity, we do not visualize flow imbalances for transitions due to teleportation.}
(B) Their homology changes with $\alpha$, which is summarized by persistence barcodes.
}
\label{fig:homo}
\end{figure*}

\begin{figure}[b]
\centering
\includegraphics[width=1\linewidth]{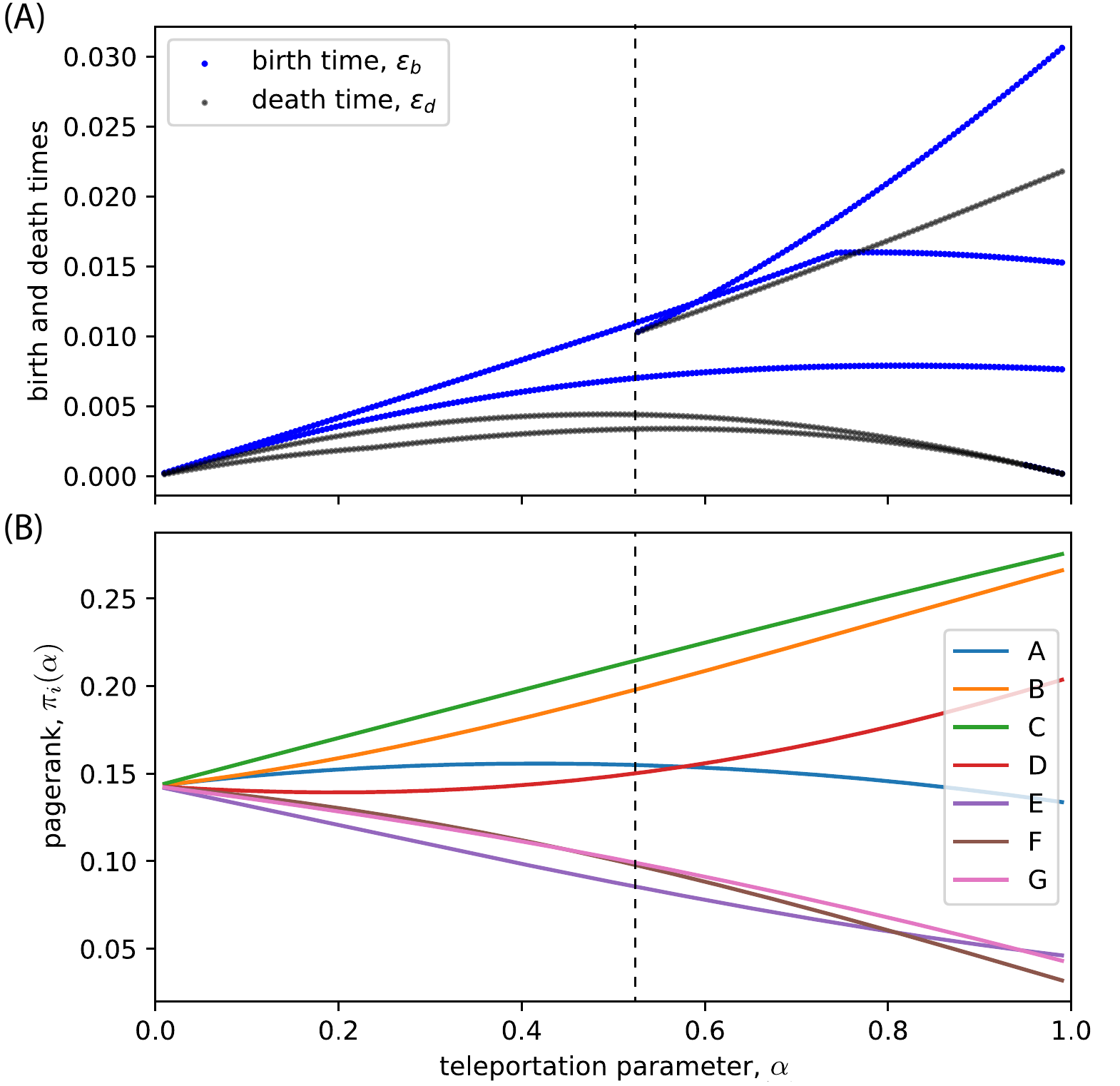}
\caption{
{\bf Bifurcation diagram summarizes homological changes onset by $\alpha$.} 
(A) Birth and death times  $\epsilon_{b,d}$ of 1-cycles arising under PageRank versus $\alpha$.  
(B) For comparison, we depict the vertices' PageRanks $\pi_i(\alpha)$.
Vertical dashed lines near $\alpha^* =0.54$ highlight that there are three cycles when $\alpha>\alpha^*$, but only two when $\alpha<\alpha^*$.
}
\label{fig:pageper}
\end{figure}

\subsection{Teleportation is a homology regularizer for PageRank}\label{sec:r1_b}

We now  study the persistent homology of convection cycles arising for the PageRank algorithm \cite{page1999pagerank,langville2006updating}, which is a popular technique to rank the importance of vertices in graphs. It has been applied to numerous applications (see survey \cite{gleich2015pagerank}), but most notably, for many years it was utilized by Google to rank website and facilitate web search.
The PageRank of a vertex $i$ is given by the stationary density $\pi_i(\alpha)$ of MC having a transition matrix of the form
\begin{align}\label{eq:pagerank}
{\bf P}(\alpha) &\equiv \alpha {\bf P} + (1-\alpha)N^{-1}{\bf 1}{\bf 1}^T,
\end{align}
where $ {\bf P}$ is the transition matrix described in Sec.~\ref{sec:conv_cycle} and $\alpha\in(0,1)$ is the \emph{teleportation parameter}. As $\alpha\to 1$, the second term vanishes and ${\bf P}(\alpha)\to{\bf P}$. Usually, $\alpha $ is chosen to be near 1 (often 0.85) so that the second term can be considered as a small perturbation that improves the mathematical characteristics of ${\bf P}$---or more formally, it is a ```regularization'' of matrix ${\bf P}$. In particular, when $\alpha \in (0,1)$ the matrix ${\bf P}(\alpha)$ is guaranteed to be  irreducible, aperiodic and positive, and  the Perron-Frobenius theorem \cite{bapat1998max} ensures that its dominant left eigenvector $\pi$ is unique and has positive entries (i.e., $ \pi_i(\alpha)>0$ for all $i$). In other words, the PageRanks are well-defined for all vertices.







We now show that the introduction of teleportation also regularizes the homology of convection cycles. In this experiment, we construct EVC filtrations with the filtration function  $f(i,j) = |\Delta_{ij}(\alpha)|$, which now depends on $\alpha$.
%
In Fig.~\ref{fig:homo}(A), we illustrate for several choices of $\alpha$  the flow imbalances that arise under PageRank, which we apply to the graph from Fig.~\ref{fig:direct}(A).
%
%
In Fig.~\ref{fig:homo}(B), we visualize their associated persistence barcodes, which we create using EVC filtrations.
Note that the choice $\alpha=1$ recovers the transition matrix, stationary distribution, flow imbalances, and persistence barcodes that were were previously studied in Figs.~\ref{fig:direct} and \ref{fig:Pbnet}.


Observe that the homological patterns of convection cycles significantly change with $\alpha$. For example, when $\alpha$ is sufficiently small, the \drt{homological}
1-cycle $\{A,B,C,D\}$ vanishes---it is ``washed out'' by the introduction of teleportation. In other word, $\alpha$ is a \emph{homology regularizer}.
This is further illustrated in Fig.~\ref{fig:pageper}(A), where we plot the birth and death times of \drt{homological} 1-cycles versus $\alpha$. For comparison, we also plot the PageRanks $\pi_i(\alpha)$ in Fig.~\ref{fig:pageper}(B).
The vertical line highlights that one of the 1-cycles vanishes when $\alpha$ decreases (approximately) below  $\alpha^* =0.54$


\drt{In Appendix~\ref{sec:appendix}, we present additional experiments that explore convection cycles arising under PageRank with $\alpha=0.8$. We show that homological 1-cycles arising for EVC filtrations with decreasing filtration parameter $\epsilon$ reveal patterns for large-flow convection cycles. In contrast, when EVC filtrations are constructed with increasing $\epsilon$, we find that the resulting homological 1-cycles     relate to small-flow convection cycles, and in particular, those involving low probability teleportation transitions.
}

\begin{figure}[t]
\centering
\includegraphics[width=1\linewidth]{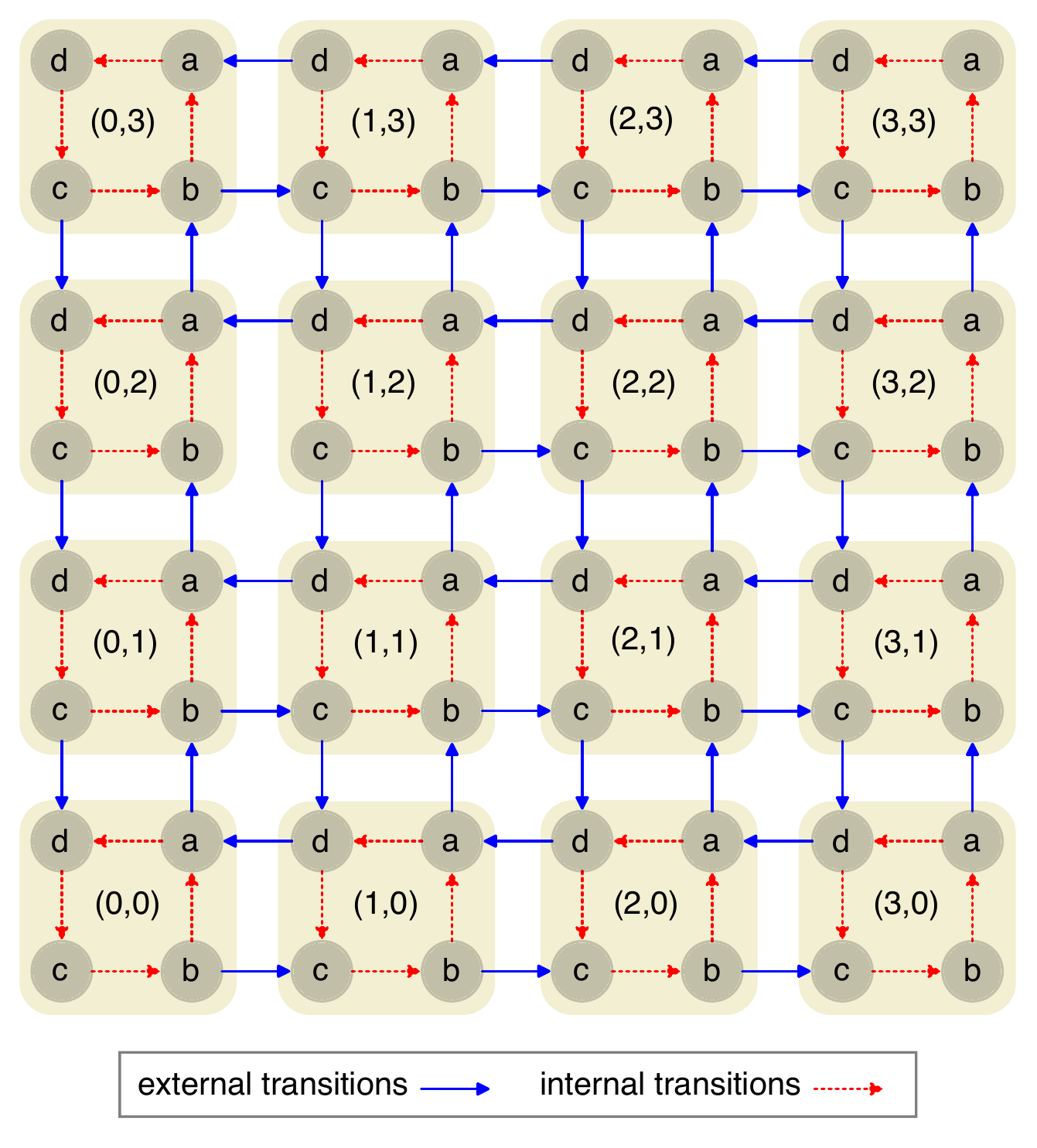}
\caption{
{\bf MC model for configuration dynamics of a  monomer.} 
The system is contains two monomers of sizes $s_1$ and $s_2$, respectively,  giving the external state $(s_1,s_2)$. Moreover, there are four internal states: a, b, c, and d. Transitions involving changes to external and internal states occur at rates $\gamma_{ex}$ and $\gamma_{in} $, respectively. 
}
\label{fig:full_lattice}
\end{figure}

\begin{figure}[t]
\centering
\includegraphics[width=1\linewidth]{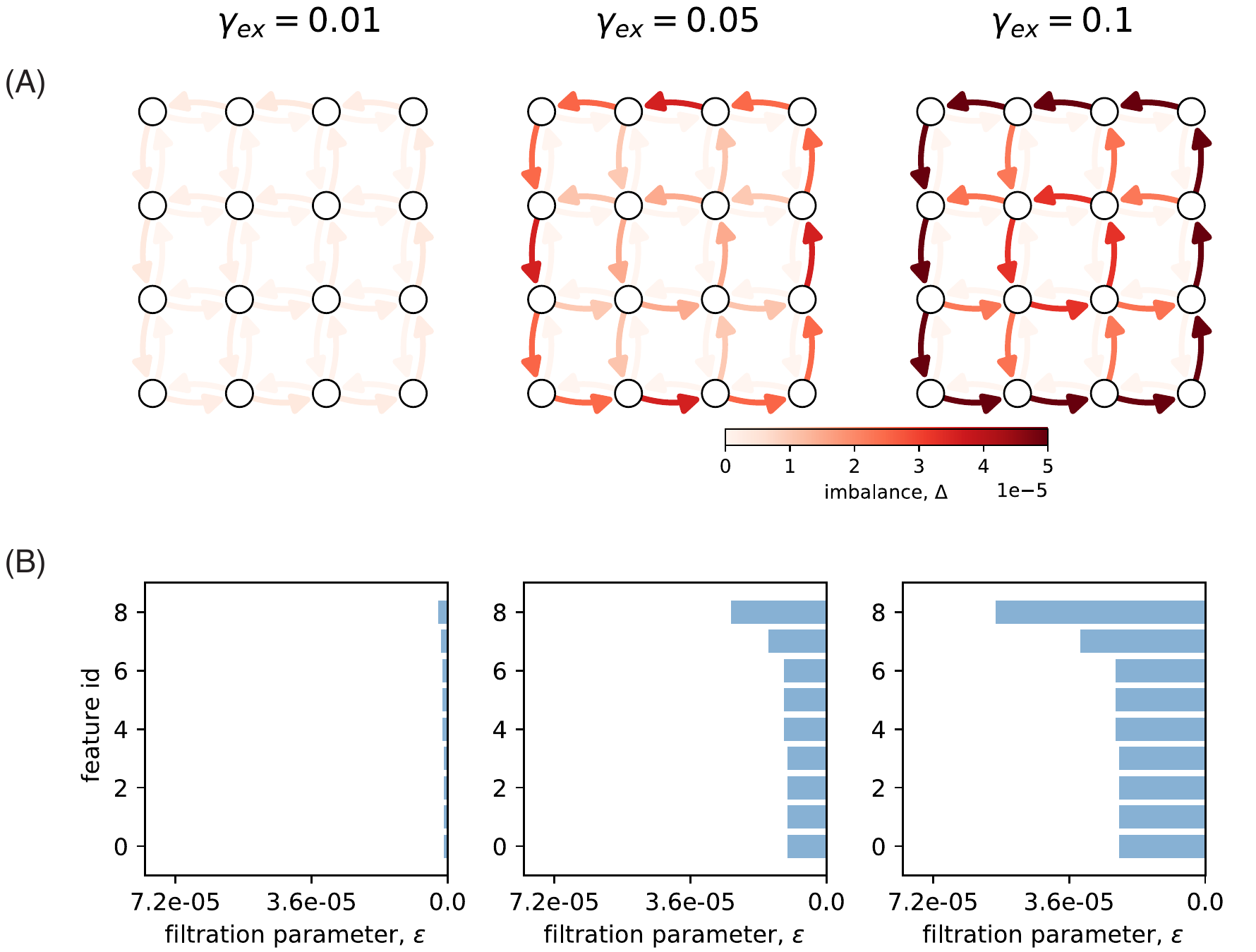}
\caption{
{\bf Persistent homology  chiral edge flow.} 
(A) Flow imbalances $\Delta_{ij}$ between external states for the bi-monomer shown in Fig.~\ref{fig:full_lattice} with  $\gamma_{in}=0.01$ and different $\gamma_{ex}$. In the limit $\gamma_{ex}\gg\gamma_{in}$ \cite{tang2021topology}, there is an emergence of a chiral edge flow, i.e., a convection cycle around the lattice's outer boundary.
(B) The corresponding persistence barcodes capture the emergence of this prominent convection cycle and other  convection cycles within the lattice.
}
\label{fig:cycle_chiral}
\end{figure}

\subsection{Persistent homology of chiral edge flows}\label{sec:r1_c}

\drt{Our second application investigates}
homological patterns of convection cycles that arise for an MC that models the stochastic configuration dynamics of two monomers.  We adopt the same notation as in \cite{tang2021topology}, which motivated our experiment. The monomer configuration (i.e., ``external state'')  is given by the number of monomers of each type, $(s_1,s_2)$, whereas the ``internal state'' is one of four possibilities: a, b, c, or d. Transitions that involve a change of internal state occur at rate $\gamma_{in}$, whereas transitions between involving external states (i.e., the addition or removal of a monomer) occur at rate $\gamma_{ex}$. The resulting MC can be visualized as a 2-dimensional lattice, which we visualize in Fig.~\ref{fig:full_lattice}. 

In Fig.~\ref{fig:cycle_chiral}(A), we visualize flow imbalances for transitions between the external states. We fix $\gamma_{in} =0.01$ and consider several $\gamma_{ex}$.  Observe that as $\gamma_{ex}$ increases, a large counter-clockwise convection cycle emerges on the boundary (i.e., ``edge'') of the lattice. This type of convection cycle is called a \emph{chiral edge flow}, and such convection cycles have important implications for  the quantum Hall effect, biological rhythms, and the dynamics of monomers \cite{tang2021topology}.
%
In Fig.~\ref{fig:cycle_chiral}(B), we visualize persistence barcodes for EVCs filtrations constructed using the method that we described in Sec.~\ref{sec:r1_a}. The chiral edge cycle corresponds to the 1-cycle with having the largest lifespan, and its homology becomes more persistent (i.e., prominent) in the limit $\gamma_{ex}\gg \gamma_{in}$.

\section{Discussion}\label{sec:conclusion}
\drt{In this paper, we examined} the patterns of convection cycles that arise under irreversible Markov chains \drt{(MCs)} from the perspective of persistent homology. Our approach required formalizing a type of filtration (EVC filtration) for scalar functions that are defined on the edges of a graph, and we studied convection cycles by choosing the filtration function to be \drt{an} MC's  flow imbalances in the stationary state. Because  Markov chains are crucial to so many diverse applications, we expect our methods to be broadly applicable across the sciences and engineering. Herein, we highlighted two such applications: the PageRank algorithm for centrality analysis and   chiral edge flows that  arise for the configuration dynamics of monomers. 
Our experiments revealed how system properties can act as homology regularizers of 
convection cycles, and we introduced homological bifurcation diagrams to summarize these changes.
This approach automates the detection, summary, and examination of convection cycles over networks, places it on stronger mathematical and computational foundations, and paves the way for further investigation into convective flows on networks.



\drt{
Additionally, our work highlights the need for new persistent homology methods to study convection cycles as well as other functions and signals defined on directed graphs. In Sec.~\ref{sec:interp}, we discuss the relation between convection cycles and homological 1-cycles, and we showed that these are two closely related, but notably different, notions of cycles. 
Sometimes there is a one-to-one correspondence between these cycles, and sometimes the relation is more complicated, due in part to the fact that a given homological k-cycle can be equivalently represented by possibly more than one homological generator. Such generators may or may not correspond to convection cycles. Moreover, convection cycles can also correspond to the boundaries of 2-simplices, and as such, they  will not be identified via the traditional tools of persistent homology.
%
%
%
Developing  persistent homology techniques that cater to convection cycles, and which specifically account for edge directions, remains an important open challenge for the applied mathematics and physics communities.
}

Our work opens up several \drt{other} new lines of research that are \drt{also}  worth noting. First, convection cycles were recently found to be an emergent property of multiplex Markov chains \cite{taylor2020multiplex} in which a set of (intralayer) Markov chains are coupled together by another set of (interlayer) Markov chains. It would be interesting to employ persistent homology  to gain a deeper understanding of this phenomenon. Second, chiral edge flows are known to be important to other applications including the quantum Hall effect and biological rhythms \cite{tang2021topology}, and future work could utilize our methods to investigate these exciting applications. Notably, our methods can reveal convection cycles that exists in addition to a chiral edge flow, which may lead to new insights for these  applications and other applications (e.g., reinforcement learning) that rely on irreversible Markov chains.

See \cite{minh_code} for a codebase that reproduces our results and can be used to study the persistent homology for convection cycles arising for other applications.

\begin{figure*}[t]
\centering
\includegraphics[width=.99\linewidth]{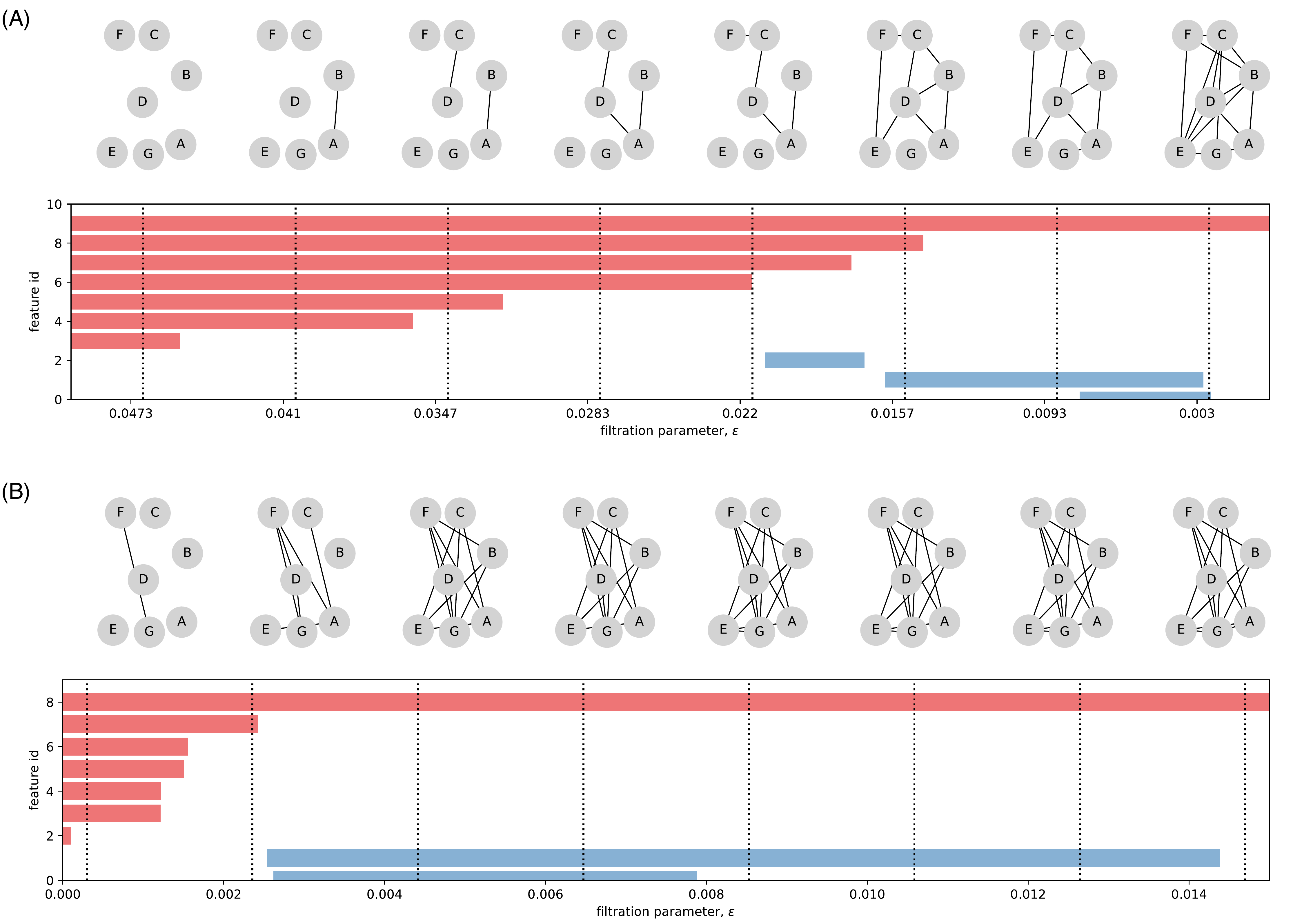}
\caption{\drt{{\bf Comparing EVC filtrations with decreasing and increasing filtration parameter $\epsilon$.}
Extending our study   in Sec.~\ref{sec:r1_b} that uses persistent homology to study convection cycles arising for a MC under the PageRank algorithm with $\alpha=0.8$, we now study homological $1$-cycles obtained via two different EVC filrations.
(A) Similar to our results in Fig.~\ref{fig:homo}, we construct EVC filtrations by  including edges for which $|\Delta_{ij}|>\epsilon$ while decreasing    $\epsilon$. Observe that the 1-cycles reveal large-flow convection cycles that are associated with large values of $|\Delta_{ij}|$.
%
(B) For comparison, we construct EVC filtrations by including weighted edges $|\Delta_{ij}|<\epsilon$ while increasing  $\epsilon$.
Observe that these 1-cycles now correspond to  small-flow convection cycles that are associated with small values of $|\Delta_{ij}|$. They primarily describe low-probability transitions that occur due to teleportation. In this work we focus on EVC filtrations with decreasing $\epsilon$, since we consider high-flow convection cycles to be the ones that are most important.
}}
\label{fig:Pbnet2}
\end{figure*}

\begin{appendix}

\section{Convection cycles are better revealed by filtrations that decrease the filtration parameter  $\epsilon$ versus increase  $\epsilon$}\label{sec:appendix}

\drt{
In Sec.~\ref{sec:PHw}, we defined EVC filtrations in which one decreases a filtration parameter $\epsilon$, retaining edges for which $f(i,j)>\epsilon$. 
Our numerical experiments that study convection cycles using persistent homology use this approach and let the filtration  be given by the   flow imbalances $f(i,j)=|\Delta_{ij}|$.
By decreasing $\epsilon$, the cycles that are first revealed correspond to large-flow convection cycles, which we consider to be the ones that are more significant.
One could also construct EVC filtrations by increasing $\epsilon$ and retaining edges for which $f(i,j)<\epsilon$. Here, we show that this latter filtration reveals  1-cycles that  relate to small-flow convection cycles, which we consider to be insignificant.

In Fig.~\ref{fig:Pbnet2}, we study EVC filtrations applied to flow imbalances arising under the PageRank algorithm with $\alpha=0.8$ for the same MC that we investigated in Sec.~\ref{sec:r1_b}. 
In Fig.~\ref{fig:Pbnet2}(A) and Fig.~\ref{fig:Pbnet2}(B), we illustrate EVC filtrations with decreasing and increasing $\epsilon$, respectively.
Observe in   Fig.~\ref{fig:Pbnet2}(A) that the  1-cycles revealed by decreasing $\epsilon$  correspond to large-flow convection cycles. In contract, observe   in Fig.~\ref{fig:Pbnet2}(B) that the  1-cycles revealed by increasing $\epsilon$ are small-flow cycles that relate to low-probability transitions that occur due to teleportation.
}


\end{appendix}

\begin{acknowledgements}
MQL and DT were supported in part by the National Science Foundation (DMS-2052720 and EDT-1551069) and the Simons Foundation (grant \#578333).
\end{acknowledgements}

\bibliographystyle{plain}
\bibliography{reference}

\begin{thebibliography}{10}

\bibitem{bapat1998max}
Ravindra~B Bapat.
\newblock A max version of the perron-frobenius theorem.
\newblock {\em Linear Algebra and its Applications}, 275:3--18, 1998.

\bibitem{edelsbrunner2010computational}
Herbert Edelsbrunner and John Harer.
\newblock {\em Computational topology: an introduction}.
\newblock American Mathematical Soc., 2010.

\bibitem{gabella2019topology}
Maxime Gabella.
\newblock Topology of learning in artificial neural networks.
\newblock {\em arXiv preprint arXiv:1902.08160}, 2019.

\bibitem{gilks1995introducing}
Walter~R Gilks, Sylvia Richardson, and David~J Spiegelhalter.
\newblock Introducing markov chain monte.
\newblock {\em Markov chain Monte Carlo in practice}, page~1, 1995.

\bibitem{giusti2015clique}
Chad Giusti, Eva Pastalkova, Carina Curto, and Vladimir Itskov.
\newblock Clique topology reveals intrinsic geometric structure in neural
  correlations.
\newblock {\em Proceedings of the National Academy of Sciences},
  112(44):13455--13460, 2015.

\bibitem{gleich2015pagerank}
David~F Gleich.
\newblock Pagerank beyond the web.
\newblock {\em siam REVIEW}, 57(3):321--363, 2015.

\bibitem{ichinomiya2020protein}
Takashi Ichinomiya, Ippei Obayashi, and Yasuaki Hiraoka.
\newblock Protein-folding analysis using features obtained by persistent
  homology.
\newblock {\em Biophysical Journal}, 118(12):2926--2937, 2020.

\bibitem{kaczynski2006computational}
Tomasz Kaczynski, Konstantin Mischaikow, and Marian Mrozek.
\newblock {\em Computational homology}, volume 157.
\newblock Springer Science \& Business Media, 2006.

\bibitem{kendall1953stochastic}
David~G Kendall.
\newblock Stochastic processes occurring in the theory of queues and their
  analysis by the method of the imbedded markov chain.
\newblock {\em The Annals of Mathematical Statistics}, pages 338--354, 1953.

\bibitem{khasawneh2016chatter}
Firas~A Khasawneh and Elizabeth Munch.
\newblock Chatter detection in turning using persistent homology.
\newblock {\em Mechanical Systems and Signal Processing}, 70:527--541, 2016.

\bibitem{kilic2022simplicial}
Bengier~Ulgen Kilic and Dane Taylor.
\newblock Simplicial cascades are orchestrated by the multidimensional geometry
  of neuronal complexes.
\newblock {\em arXiv preprint arXiv:2201.02071}, 2022.

\bibitem{kingman1969markov}
John~FC Kingman.
\newblock Markov population processes.
\newblock {\em Journal of Applied Probability}, pages 1--18, 1969.

\bibitem{kondic2012topology}
L~Kondic, A~Goullet, CS~O'Hern, M~Kramar, Konstantin Mischaikow, and
  RP~Behringer.
\newblock Topology of force networks in compressed granular media.
\newblock {\em EPL (Europhysics Letters)}, 97(5):54001, 2012.

\bibitem{kramar2013persistence}
M~Kramar, Arnaud Goullet, Lou Kondic, and Konstantin Mischaikow.
\newblock Persistence of force networks in compressed granular media.
\newblock {\em Physical Review E}, 87(4):042207, 2013.

\bibitem{langville2006updating}
Amy~N Langville and Carl~D Meyer.
\newblock Updating markov chains with an eye on google's pagerank.
\newblock {\em SIAM journal on matrix analysis and applications},
  27(4):968--987, 2006.

\bibitem{minh_code}
Minh~Quang Le.
\newblock Codebase for persistent homology of convection cycles in network
  flows
  \url{https://github.com/minhquan89/Persistent-Homology-of-Convection-Cycles}.

\bibitem{liang1998analytical}
Jie Liang, Herbert Edelsbrunner, Ping Fu, Pamidighantam~V Sudhakar, and Shankar
  Subramaniam.
\newblock Analytical shape computation of macromolecules: I. molecular area and
  volume through alpha shape.
\newblock {\em Proteins: Structure, Function, and Bioinformatics}, 33(1):1--17,
  1998.

\bibitem{liu2019scalable}
Shusen Liu, Di~Wang, Dan Maljovec, Rushil Anirudh, Jayaraman~J Thiagarajan,
  Sam~Ade Jacobs, Brian~C Van~Essen, David Hysom, Jae-Seung Yeom, Jim Gaffney,
  et~al.
\newblock Scalable topological data analysis and visualization for evaluating
  data-driven models in scientific applications.
\newblock {\em IEEE transactions on visualization and computer graphics},
  26(1):291--300, 2019.

\bibitem{lovasz1993random}
L{\'a}szl{\'o} Lov{\'a}sz et~al.
\newblock Random walks on graphs: A survey.
\newblock {\em Combinatorics, Paul erdos is eighty}, 2(1):1--46, 1993.

\bibitem{motta2019hyperparameter}
Francis Motta, Christopher Tralie, Rossella Bedini, Fabiano Bini, Gilberto
  Bini, Hamed Eramian, Marcio Gameiro, Steve Haase, Hugh Haddox, John Harer,
  et~al.
\newblock Hyperparameter optimization of topological features for machine
  learning applications.
\newblock In {\em 2019 18th IEEE International Conference On Machine Learning
  And Applications (ICMLA)}, pages 1107--1114. IEEE, 2019.

\bibitem{Note1}
The notion of ``dimension'' is more interpretable for the case of a
  non-abstract simplicial complex, for which the vertices correspond to
  locations in a Euclidean metric space. In that case, every $k$-simplex is
  defined as the $k$-dimensional surface that is contained by its faces, which
  themselves are $(k-1)$ dimensional surfaces. For example, a 2-simplex is a
  triangle defined as the interior of 3 line segments (which are the 3
  1-dimensional cofaces of the 2-simplex).

\bibitem{otter2017roadmap}
Nina Otter, Mason~A Porter, Ulrike Tillmann, Peter Grindrod, and Heather~A
  Harrington.
\newblock A roadmap for the computation of persistent homology.
\newblock {\em EPJ Data Science}, 6(1):17, 2017.

\bibitem{page1999pagerank}
Lawrence Page, Sergey Brin, Rajeev Motwani, and Terry Winograd.
\newblock The pagerank citation ranking: Bringing order to the web.
\newblock Technical report, Stanford InfoLab, 1999.

\bibitem{parr1998reinforcement}
Ronald Parr and Stuart Russell.
\newblock Reinforcement learning with hierarchies of machines.
\newblock {\em Advances in neural information processing systems}, pages
  1043--1049, 1998.

\bibitem{perea2015sliding}
Jose~A Perea and John Harer.
\newblock Sliding windows and persistence: An application of topological
  methods to signal analysis.
\newblock {\em Foundations of Computational Mathematics}, 15(3):799--838, 2015.

\bibitem{petri2014homological}
Giovanni Petri, Paul Expert, Federico Turkheimer, Robin Carhart-Harris, David
  Nutt, Peter~J Hellyer, and Francesco Vaccarino.
\newblock Homological scaffolds of brain functional networks.
\newblock {\em Journal of The Royal Society Interface}, 11(101):20140873, 2014.

\bibitem{tang2021topology}
Evelyn Tang, Jaime Agudo-Canalejo, and Ramin Golestanian.
\newblock Topology protects chiral edge currents in stochastic systems.
\newblock {\em Physical Review X}, 11(3):031015, 2021.

\bibitem{taylor2020multiplex}
Dane Taylor.
\newblock Multiplex markov chains: Convection cycles and optimality.
\newblock {\em Physical Review Research}, 2(3):033164, 2020.

\bibitem{taylor2015topological}
Dane Taylor, Florian Klimm, Heather~A Harrington, Miroslav Kram{\'a}r,
  Konstantin Mischaikow, Mason~A Porter, and Peter~J Mucha.
\newblock Topological data analysis of contagion maps for examining spreading
  processes on networks.
\newblock {\em Nature communications}, 6:7723, 2015.

\bibitem{gudhi:urm}
{The GUDHI Project}.
\newblock {\em {GUDHI} User and Reference Manual}.
\newblock {GUDHI Editorial Board}, {3.4.1} edition, 2021.

\bibitem{tierney1994markov}
Luke Tierney.
\newblock Markov chains for exploring posterior distributions.
\newblock {\em the Annals of Statistics}, pages 1701--1728, 1994.

\bibitem{xin1994topology}
Ying-Jie Xin and Yuan-Hua Zhou.
\newblock Topology on image processing.
\newblock In {\em Proceedings of ICSIPNN'94. International Conference on
  Speech, Image Processing and Neural Networks}, pages 764--767. IEEE, 1994.

\end{thebibliography}

\end{document}